1# Simultaneous Robust State Estimation, Topology Error Processing, and Outage Detection for Unbalanced Distribution Systems

Zahra Soltani, *Student member, IEEE*, Shanshan Ma, *Member, IEEE*, Mojdeh Khorsand, *Member, IEEE*, and Vijay Vittal, *Fellow, IEEE*

*Abstract*—This paper proposes an efficient algorithmic approach that overcomes the critical challenges in the real-time unbalanced distribution system state estimation, topology error processing, and outage identification simultaneously: (1) Limited locations of measurement devices and unsynchronized measurement data as well as missing and bad data, (2) Complicated mixed-phase switch actions and mutual impedances, and (3) the nonlinear nature of unbalanced distribution system power flow and measurement data . A single snap-shot mixed-integer quadratic programming (MIQP) optimization framework is proposed to cope with these challenges. This MIQP framework presents a more accurate unbalanced distribution system model, linearizes the nonlinear relationship between bus voltage and current injections, formulates the complicated mixed-phase switch operations, and executes the outage detection function via analytic constraints. The proposed model can be effectively solved by a general commercial MIP solver. The effectiveness of the proposed approach is verified on an actual distribution feeder in Arizona.

*Index Terms*— Distribution system state estimation, topology error process, outage detection, unbalanced distribution systems, AC optimal power flow, mutual impedances.

NOMENCLATURE

*Sets and Indices*

| | |
|---|---|
| $\mathcal{J}/n, m$ | Set/Index for buses |
| $\mathcal{H}/\mathcal{L}$ | Set/Index of lines |
| $\mathcal{H}'/\mathcal{H}''$ | Set of lines with/without switch |
| $\psi/\phi$ | Set/Index of phases |
| $\psi/\phi$ | Set/Index of phases |

*Parameters and Constants*

| | |
|---|---|
| $Z_{\mathcal{L}}^{\phi,p}/y_{\mathcal{L}}^{\phi,p}$ | Impedance/ Shunt admittance of line $\mathcal{L}$ between phases $\phi$ and $p$ |
| $R_{\mathcal{L}}^{\phi,p}/X_{\mathcal{L}}^{\phi,p}$ | Resistance/Reactance of line $\mathcal{L}$ between phases $\phi$ and $p$ |
| $NL_{r,\phi}^{TR}$ | No load loss transformer $r$ at phase $\phi$ |
| $M$ | A large positive number |

*Variables*

| | |
|---|---|
| $I_{\mathcal{L}}^{\phi}$ | Current of phase $\phi$ of line $\mathcal{L}$ |
| $I_{\mathcal{L}}^{r,p}/I_{\mathcal{L}}^{im,p}$ | Real/Imaginary part of phase $\phi$ current of line $\mathcal{L}$ |
| $V_n^{\phi}$ | Voltage at bus $n$ and phase $\phi$ |
| $V_n^{r,\phi}/V_n^{im,\phi}$ | Real/ Imaginary part of voltage at bus $n$ and phase $\phi$ |
| $I_n^{r,\phi}/I_n^{im,\phi}$ | Real/Imaginary part of current injection at bus $n$ and phase $\phi$ |
| $P_{g,\phi}^G/Q_{g,\phi}^G$ | Active/Reactive power of DER $g$ at phase $\phi$ |
| $d_{p,l,\phi}$ | Active power of load $l$ at phase $\phi$ |
| $d_{q,l,\phi}$ | Reactive power of load $l$ at phase $\phi$ |
| $Q_{c,\phi}^C$ | Reactive power of capacitor unit $c$ at phase $\phi$ |
| $u_{\mathcal{L}}^{\phi}$ | Status of switch of line $\mathcal{L}$ at phase $\phi$ |
| $\Lambda_{\mathcal{L}}^k$ | Auxiliary binary variable |

## I. Introduction

THE variabilities and uncertainties of the high penetration of distributed energy resources (DERs) and demand-respond devices have further complicated the electric distribution system operation, especially the network topology changes frequently due to frequent reverse flow and protection actions [1], [2]. To monitor distribution system operation under various variabilities and uncertainties, a real-time accurate distribution system state estimation (DSSE) is the crucial component of the distribution management system (DMS) to provide inputs for voltage and power flow control [1], [3], [4]. However, the current DSSE is based on a strong assumption that DMS has precise topology knowledge [3], which indicates that the topology processor has to provide correct topology identification information in advance. This strong assumption may leave DSSE vulnerable to topology identification errors [5]. Since most switches are still operated manually by field crews and have no communication links to automatically send the switch status back to DMS. At the same time, the complicated mixed-phase multiple switch actions are easily triggered frequently to satisfy DERs or demand response operation and protection requirements, which may not be detected and generate topology errors for DSSE. Furthermore, the limited deployment number of micro-Phasor Measurement Units (micro-PMUs) in the current distribution system and the unsynchronized measurement data from smart meter data at the customer locations cannot provide real-time complete observability information for DMS. Therefore, it is necessary to develop an efficient real-time distribution system topology processor and state estimation tool to correctly identify the

The work is funded by the Department of Energy (DOE) Advanced Research Projects Agency–Energy (ARPA-E) under OPEN 2018 program. Zahra Soltani, Shanshan Ma, Mojdeh Khorsand, and Vijay Vittal are with the School of Electrical, Computer, and Energy Engineering, Arizona State University, Tempe, AZ 85281 USA.

network topology configuration and estimate the system state with limited observation information.

Many efforts have been made to investigate studies on the distribution system topology processor. In [6], the authors presented a time-series signature verification method to detect distribution network topology. It explores possible topology configurations of a given distribution network to build a signal library and projects actual voltage phasorial patterns onto the corresponding library. However, this topology verification method depends on prior information of switch status and voltage measurement from micro-PMUs and three-parameter tuning. If some errors in the prior switch data or the variation of DER output and load become larger, it may generate errors in the current topology identification results. A statistical inference framework was proposed in [7] to verify the single-phase distribution grid topology by exploiting the voltage covariance from smart meter data. It fails to consider the unbalanced switch actions or multi-phase switch actions in the distribution system. A deep-learning-based topology identification approach was developed in [8] to identify multiple multi-phase switch actions with high accuracy. However, this deep-learning-based method requires massive data to train its identifier model. In summary, [5-8] identified topology under normal conditions but did not explore the possibility of outage detection in their researches. Although a MILP formulation was developed in [9] to identify switch malfunctions and detect outages, this formulation is still based on balanced power flow, which also fails to consider the unbalanced mixed-phase switch operations in the unbalanced distribution systems.

In addition to distribution system state estimation (DSSE), the latest approaches include optimization-based and data-driven methods. A weighted least square (WLS) based DSSE model was proposed in [10] for three-phase four-conductor configured unsymmetrical medium voltage distribution systems. Although [10] developed an efficient state variable reduction approach to improve DSSE computational performance, the nonlinear nature in the WSL-based DSSE problem and unbalanced three-phase power flow makes this problem hard to obtain the global optimal solution, even provides an acceptable solution for a large-scale distribution system. The authors in [11] developed a linear DSSE algorithm using the Taylor series of voltages in the interval form and solved by interval arithmetic techniques. However, this linear DSSE algorithm does not consider the unbalanced multi-phase power flow feature of distribution systems. It is also not appropriate to assume the node voltage magnitude is approximate to be the real part of voltage concerning the mixed-phase voltage. In [12], the nonconvex DSSE problem was reformulated as a rank-constrained Semidefinite programming problem and can be solved by the rank reduction approach and the convex iteration approach. [13] presented a data-driven learning-based optimization approach for the WLS-based DSSE problem. The proposed approach used historical load and energy generation data to train a shallow neural network to learn good initialization points for the Gauss-Newton algorithm. A deep learning approach was proposed in [14] to Bayesian state estimation for an unobservable real-time distribution system. However, data-driven methods in [13], [14] rely on historical data and may not adapt fast to sudden network changes.

As the observability of distribution systems has been increased by using uPMUs, Intelligent Electronic Devices (IEDs), voltage regulators, and smart inverters of DER [15], some researchers integrate the topology processor with distribution system state estimation together. An algorithm based on changing the branch statuses one after the other and performing branch current state estimation in each case was proposed in [16] to identify distribution system topology errors. Authors in [17] parallelly run distribution system state estimation model with possible critical network configuration using real measurement, use a recursive Bayesian approach to build a probability model, and identify distribution network configuration by the highest probability. However, before algorithms in [16, 17] can identify topology errors, both of them need to run multiple state estimation cases first, which may not be efficient when multiple switch actions in different phases happen simultaneously. In [18], a generalized state estimation approach was extended to identify distribution system topology changes by formulating the switch statuses as continuous state variables and adding additional soft operational constraints for each switch device. But this formulation also does not consider the unbalanced multi-phase switch actions in the unbalanced distribution systems.

To deal with the deficiencies of these methods mentioned above, this paper proposes a single snap-shot MIQP formulation to identify distribution system switch actions, estimate system state, and detect outages simultaneously. To accurately model the unbalanced distribution systems, a current and voltage (IV) based AC optimal power flow (ACOPF) considering the mutual impedance and admittance among different phases is proposed. The proposed ACOPF based on IV format is novel and different from the widely-applied unbalanced DistFlow model, which assumes that unbalanced three-phase voltages are nearly balanced and ignore the power losses caused by high X/R ratio in [19]. An iterative first-order approximation of the Tylor extension approach is inserted to linearize the nonlinear power balance injection constraints in the proposed IV-based ACOPF model. At the same time, the complicated mixed-phase switch actions are formulated using the big-M method. The outage detection function is created by analytic constraints. An unsynchronized data integration approach is developed to deal with the smart meter and uPMU data unsynchronized issues. Simultaneously, the issue that smart meter data does not contain reactive power information is solved by adding a power factor constraint in the formulation. The robustness of the proposed tool is evaluated against bad data in the sensors' measughrments. The main contribution of this paper can be summarized as follows:

(i) The proposed detailed unbalanced distribution system modeling based on IV-based ACOPF can provide more accurate DSSE information for DMS, especially the voltage profile.
(ii) Different techniques are proposed to linearize the mixed-integer nonlinear joint DSSE and topology processor

problem into a MIQP optimization problem. The peoposed MIQP model can be solved by a commercial solver (e.g., Gurobi and CPLEX), which avoid numerical instability and initial point sensitivity caused by the nonlinearity.
(iii) The MIQP formulation can avoid multiple runs of the off-line WLS-based DSSE procedure and do not rely on prior or historical information.
(iv) The analytic formulation of the outage detection function and mixed-phase switch operation can help DMS to detect dynamic switch operation or fault accurately in real-time.
(v) The proposed unsynchronized data integration approach can make full use of the real-time data.

The paper has the following structure. The mathematic formulation of the simultaneous topology processor and state estimation tool is modeled in Section II. The analytic formulation of outage detection inserted in the proposed tool is developed in Section III. Section IV presents an unsynchronous sensor data integration approach. Simulation results are illustrated in Section V. Concluding remarks are summarized in Section VI.

## II. SIMULTANEOUS TOPOLOGY PROCESSOR AND STATE ESTIMATION TOOL IN UNBALANCED DISTRIBUTION NETWORKS

In this section, a new ACOPF based on IV formulation for the multi-phase unbalanced distribution system is proposed, which is utilized in the simultaneous topology processor and state estimation tool. Then, the nonlinear and proposed models of the simultaneous topology processor and state estimation tool, which can be executed in the DMS, are presented.

### A. The Unbalanced Distribution System ACOPF Modeling

Let consider an unbalanced distribution network with graph representation of $\mathcal{G} = (\mathcal{J}, \mathcal{H})$. $\mathcal{J}$ is a set of buses and $\mathcal{H}$ is a set of distribution lines. Every bus and line in $\mathcal{G}$ can have three phases according to set $\psi = \{a, b, c\}$. $\mathcal{H}$ includes a set of lines without switch (denoted by $\mathcal{H}'$) and a set of lines with switch (denoted by) $\mathcal{H}''$. In the proposed formulation, inverse of $x$ is represented by $(x)^{-1}$. The voltage difference over a three-phase distribution line $\mathcal{L} \in \mathcal{H}$ at phase $\phi$ is formulated based on self- and mutual impedances and admittances of the line as follows:

$$V_n^\phi - V_m^\phi = Z_{\mathcal{L}}^{\phi,\phi} I_{\mathcal{L}}^\phi - \frac{1}{2} Z_{\mathcal{L}}^{\phi,\phi} y_{\mathcal{L}}^{\phi,\phi} V_n^\phi -$$
$$\frac{1}{2} Z_{\mathcal{L}}^{\phi,\phi} \left( \sum_{k \in \psi, k \neq \phi} y_{\mathcal{L}}^{\phi,k} V_n^k \right) + \sum_{p \in \psi, p \neq \phi} Z_{\mathcal{L}}^{\phi,p} \left( I_{\mathcal{L}}^p - \frac{1}{2} \left( \sum_{k \in \psi} y_{\mathcal{L}}^{p,k} V_n^k \right) \right), \forall \phi \in \psi, \mathcal{L} \in \mathcal{H}' \quad (1a)$$

where line $\mathcal{L}$ connects bus $n$ to bus $m$. The current from bus $n$ to bus $m$ at the phase $\phi$ of a distribution line is defined based on (1a) as follows:

$$I_{\mathcal{L}}^\phi = \left( Z_{\mathcal{L}}^{\phi,\phi} \right)^{-1} \left[ V_n^\phi - V_m^\phi + \frac{1}{2} Z_{\mathcal{L}}^{\phi,\phi} y_{\mathcal{L}}^{\phi,\phi} V_n^\phi + \frac{1}{2} Z_{\mathcal{L}}^{\phi,\phi} \left( \sum_{k \in \psi, k \neq \phi} y_{\mathcal{L}}^{\phi,k} V_n^k \right) - \sum_{p \in \psi, p \neq \phi} Z_{\mathcal{L}}^{\phi,p} \left( I_{\mathcal{L}}^p - \frac{1}{2} \left( \sum_{k \in \psi} y_{\mathcal{L}}^{p,k} V_n^k \right) \right) \right], \forall \phi \in \psi, \mathcal{L} \in \mathcal{H}' \quad (1b)$$

As it is illustrated in (1b), the current of three-phase distribution line $\mathcal{L}$ at the phase $\phi$ is obtained not only based on voltage, impedance, and admittance values of the phase $\phi$ but also on the voltage of the other phases as a result of the mutual impedances and admittances between the phase $\phi$ and the other two phases of the line. Based on (1b), $I_{n,m}^\phi$ has five terms, where the first two terms (i.e., $\left( Z_{\mathcal{L}}^{\phi,\phi} \right)^{-1} V_n^\phi - \left( Z_{\mathcal{L}}^{\phi,\phi} \right)^{-1} V_m^\phi$) represent the current flowing in the series self-impedance (i.e., $Z_{n,m}^{\phi,\phi}$) of the three-phase distribution line. Third term (i.e., $\frac{1}{2} \left( Z_{\mathcal{L}}^{\phi,\phi} \right)^{-1} Z_{\mathcal{L}}^{\phi,\phi} y_{\mathcal{L}}^{\phi,\phi} V_n^\phi$) shows the current of self-shunt admittance (i.e., $y_{\mathcal{L}}^{\phi,\phi}$) of the line. Fourth term (i.e., $\frac{1}{2} \left( Z_{\mathcal{L}}^{\phi,\phi} \right)^{-1} Z_{\mathcal{L}}^{\phi,\phi} \left( \sum_{k \in \psi, k \neq \phi} y_{\mathcal{L}}^{\phi,k} V_n^k \right)$) shows the effect of mutual shunt admittances between phase $\phi$ and the other two phases on $I_{n,m}^\phi$. In the last term (i.e., $\left( Z_{\mathcal{L}}^{\phi,\phi} \right)^{-1} \sum_{p \in \psi, p \neq \phi} Z_{\mathcal{L}}^{\phi,p} \left( I_{\mathcal{L}}^p - \frac{1}{2} \left( \sum_{k \in \psi} y_{\mathcal{L}}^{p,k} V_n^k \right) \right)$), the current flowing in the other two phases (i.e., $I_{\mathcal{L}}^p$) and also the self- and mutual shunt admittances of the other two phases (i.e., $\frac{1}{2} \left( \sum_{k \in \psi} y_{\mathcal{L}}^{p,k} V_n^k \right)$) are appeared in the formulation of $I_{\mathcal{L}}^\phi$ due to the mutual impedances between the phase $\phi$ and the other two phases (i.e., $Z_{\mathcal{L}}^{\phi,p}$). For instance, the current at phase $a$ of a line with sending bus $n$ and receiving bus $m$ is given in (1c).

$$I_{\mathcal{L}}^a = \left( Z_{\mathcal{L}}^{a,a} \right)^{-1} \left[ V_n^a - V_m^a + \frac{1}{2} Z_{\mathcal{L}}^{a,a} y_{\mathcal{L}}^{a,a} V_n^a + \frac{1}{2} Z_{\mathcal{L}}^{a,a} \left( y_{\mathcal{L}}^{a,b} V_n^b + y_{\mathcal{L}}^{a,c} V_n^c \right) - Z_{\mathcal{L}}^{a,b} \left( I_{\mathcal{L}}^b - \frac{1}{2} \left( y_{\mathcal{L}}^{b,a} V_n^a + y_{\mathcal{L}}^{b,b} V_n^b + y_{\mathcal{L}}^{b,c} V_n^c \right) \right) - Z_{\mathcal{L}}^{a,c} \left( I_{\mathcal{L}}^c - \frac{1}{2} \left( y_{\mathcal{L}}^{c,a} V_n^a + y_{\mathcal{L}}^{c,b} V_n^b + y_{\mathcal{L}}^{c,c} V_n^c \right) \right) \right] \quad (1c)$$

The real and imaginary components of current of phase $\phi$ of line $\mathcal{L} \in \mathcal{H}'$ (i.e., $I_{\mathcal{L}}^\phi$ in (1b)) are given in (1d)-(1e), respectively.

$$I_{\mathcal{L}}^{r,\phi} = \left( R_{\mathcal{L}}^{\phi,\phi} \right)^{-1} \left[ V_n^{r,\phi} - V_m^{r,\phi} + X_{\mathcal{L}}^{\phi,\phi} I_{\mathcal{L}}^{im,\phi} - \frac{1}{2} y_{\mathcal{L}}^{\phi,\phi} \left( R_{\mathcal{L}}^{\phi,\phi} V_n^{im,\phi} + X_{\mathcal{L}}^{\phi,\phi} V_n^{r,\phi} \right) - \frac{1}{2} \left( \sum_{k \in \psi, k \neq \phi} y_{\mathcal{L}}^{\phi,k} \left( R_{\mathcal{L}}^{\phi,\phi} V_n^{im,k} + X_{\mathcal{L}}^{\phi,\phi} V_n^{r,k} \right) \right) - \sum_{p \in \psi, p \neq \phi} R_{\mathcal{L}}^{\phi,p} \left( I_{\mathcal{L}}^{r,p} + \frac{1}{2} \left( \sum_{k \in \psi} y_{\mathcal{L}}^{p,k} V_n^{im,k} \right) \right) + \sum_{p \in \psi, p \neq \phi} X_{\mathcal{L}}^{\phi,p} \left( I_{\mathcal{L}}^{im,p} - \frac{1}{2} \sum_{k \in \psi} y_{\mathcal{L}}^{p,k} V_n^{r,k} \right) \right], \forall \phi \in \psi \quad (1d)$$

$$I_{\mathcal{L}}^{im,\phi} = \left( R_{\mathcal{L}}^{\phi,\phi} \right)^{-1} \left[ V_n^{im,\phi} - V_m^{im,\phi} - X_{\mathcal{L}}^{\phi,\phi} I_{\mathcal{L}}^{r,\phi} + \frac{1}{2} y_{\mathcal{L}}^{\phi,\phi} \left( R_{\mathcal{L}}^{\phi,\phi} V_n^{r,\phi} - X_{\mathcal{L}}^{\phi,\phi} V_n^{im,\phi} \right) + \frac{1}{2} \left( \sum_{k \in \psi, k \neq \phi} y_{\mathcal{L}}^{\phi,k} \left( R_{\mathcal{L}}^{\phi,\phi} V_n^{r,k} - X_{\mathcal{L}}^{\phi,\phi} V_n^{im,k} \right) \right) - \sum_{p \in \psi, p \neq \phi} R_{\mathcal{L}}^{\phi,p} \left( I_{\mathcal{L}}^{im,p} - \frac{1}{2} \left( \sum_{k \in \psi} y_{\mathcal{L}}^{p,k} V_n^{r,k} \right) \right) - \sum_{p \in \psi, p \neq \phi} X_{\mathcal{L}}^{\phi,p} \left( I_{\mathcal{L}}^{r,p} + \frac{1}{2} \sum_{k \in \psi} y_{\mathcal{L}}^{p,k} V_n^{im,k} \right) \right], \forall \phi \in \psi \quad (1e)$$

The injected current to the phase $\phi$ of the bus $n \in \mathcal{J}$ is expressed using (1f)-(1g).

$$I_n^{r,\phi} = \sum_{m \in \delta(n)} I_{\mathcal{L}}^{r,\phi}, \forall \phi \in \psi, n \in \mathcal{J} \quad (1f)$$

$$I_n^{im,\phi} = \sum_{m \in \delta(n)} I_{\mathcal{L}}^{im,\phi}, \forall \phi \in \psi, n \in \mathcal{J} \quad (1g)$$

where $\delta(n)$ is set of buses connected to bus $n$. Let assume an unbalanced distribution system, where DERs such as solar panels are connected via distribution transformers to the network. For each phase $\phi$ of the bus $n \in \mathcal{J}$, the power balance formulations are given in (1h)-(1i).

$\sum_{\forall g \in g(n)} P_{g,\phi}^G - \sum_{\forall l \in l(n)} d_{p,l,\phi} - \sum_{\forall r \in R(n)} NL_{r,\phi}^{TR} = V_n^{r,\phi} I_n^{r,\phi} + V_n^{im,\phi} I_n^{im,\phi}, \forall \phi \in \psi, n \in \mathcal{J}$ (1h)

$\sum_{\forall g \in g(n)} Q_{g,\phi}^G + \sum_{\forall c \in C(n)} Q_{c,\phi}^C - \sum_{\forall l \in l(n)} d_{q,l,\phi} = V_i^{im,a} I_i^{r,a} - V_i^{r,a} I_i^{im,a}, \forall \phi \in \psi, n \in \mathcal{J}$ (1i)

In the most of distribution ACOPF models such as DistFlow model [19], line flow constraints are nonlinear, which are linearized by ignoring line losses in the distribution systems. In the proposed ACOPF model (1), not only line losses are modeled but also the line flow constraints are inherently linear. Therefore, the proposed ACOPF model is more accurate for the topology detection and state estimation in distribution systems with high losses.

*B. Nonlinear Simultaneous Topology Processor and State Estimation Tool*

In this section, the proposed ACOPF model in (1) is utilized for developing the simultaneous topology processor and state estimation model. For an unbalanced distribution system with a set of lines $\mathcal{H}$, where $\mathcal{H} = \mathcal{H}' \cup \mathcal{H}''$, the line current constraints for each line $\mathcal{L} \in \mathcal{H}'$ (i.e., a line without switch) are defined in (1b) and (1d)-(1e). If a switch device is installed on each phase of a three-phase distribution line, the formulation of the current at phase $\phi$ of the line presented in (1b) should be modified to model possible connectivity statuses of the switches (i.e., connected or disconnected). In order to model the status of switch devices, a binary variable is introduced for each switch installed at each phase of the three-phase line $\mathcal{L} \in \mathcal{H}''$ denoted as $u_\mathcal{L}^\phi$. In this paper, $u_\mathcal{L}^\phi = 0$ implies that the phase $\phi$ of the line $\mathcal{L}$ is disconnected, and $u_\mathcal{L}^\phi = 1$ shows that phase $\phi$ of the line $\mathcal{L}$ is connected. Therefore, the current at the phase $\phi$ of the line $\mathcal{L} \in \mathcal{H}''$ from bus $n$ to bus $m$ with three single-phase switches is defined using (6).

$I_\mathcal{L}^\phi = (Z_\mathcal{L}^{\phi,\phi})^{-1} u_\mathcal{L}^\phi \left[ V_n^\phi - V_m^\phi + \frac{1}{2} Z_\mathcal{L}^{\phi,\phi} y_\mathcal{L}^{\phi,\phi} V_n^\phi + \frac{1}{2} Z_\mathcal{L}^{\phi,\phi} (\sum_{k \in \psi, k \neq \phi} y_\mathcal{L}^{\phi,k} u_\mathcal{L}^k V_n^k) - \sum_{p \in \psi, p \neq \phi} Z_\mathcal{L}^{\phi,p} u_\mathcal{L}^p \left( I_\mathcal{L}^p - \frac{1}{2} (\sum_{k \in \psi} y_\mathcal{L}^{p,k} \Lambda_\mathcal{L}^k V_n^k) \right) \right], \forall \phi \in \psi, \mathcal{L} \in \mathcal{H}''$ (2a)

where $\Lambda_\mathcal{L}^k$ is an axillary binary variable associated with status of single-phase switches of the line as defined in (2b).

$\Lambda_\mathcal{L}^k = \begin{cases} u_\mathcal{L}^k & \text{if } k \neq \phi \vee p \\ 1 & \text{if } k = \phi \vee p \end{cases}$ (2b)

where "or" is indicated by symbol ∨. In (2b), $\Lambda_\mathcal{L}^k$ is equal to $u_\mathcal{L}^k$, when $k$ is not equal to $\phi$ (i.e., phase of $I_\mathcal{L}^\phi$) or $p$. According to (2a), $I_\mathcal{L}^\phi$ is not only determined based on the status of the switch installed at the phase $\phi$ (i.e., $u_\mathcal{L}^\phi$) but also on the status of switches installed on the other two phases because of the mutual impedances and admittances. For instance, the current at phase $a$ of the line $\mathcal{L} \in \mathcal{H}''$ with three single-phase switches is obtained using (2c).

$I_\mathcal{L}^a = (Z_\mathcal{L}^{a,a})^{-1} u_\mathcal{L}^a \left[ V_n^a - V_m^a + \frac{1}{2} Z_\mathcal{L}^{a,a} y_\mathcal{L}^{a,a} V_n^a + \frac{1}{2} Z_\mathcal{L}^{a,a} (y_\mathcal{L}^{a,b} u_\mathcal{L}^b V_n^b + y_{n,m}^{a,c} u_\mathcal{L}^c V_n^c) - Z_\mathcal{L}^{a,b} u_\mathcal{L}^b \left( I_\mathcal{L}^b - \frac{1}{2} (y_\mathcal{L}^{b,a} V_n^a + y_\mathcal{L}^{b,b} V_n^b + y_\mathcal{L}^{b,c} u_\mathcal{L}^c V_n^c) \right) - Z_\mathcal{L}^{a,c} u_\mathcal{L}^c \left( I_\mathcal{L}^c - \frac{1}{2} (y_\mathcal{L}^{c,a} V_n^a + y_\mathcal{L}^{c,b} u_\mathcal{L}^b V_n^b + y_\mathcal{L}^{c,c} V_n^c) \right) \right]$ (2c)

As shown in (2c), the current at phase $a$ of the line is not only dependent on the status of the switch of phase $a$ (i.e., $u_\mathcal{L}^a$) but also on the status of the switches of phases $b$ and $c$ (i.e., $u_\mathcal{L}^b$ and $u_\mathcal{L}^c$). Also, in order to elucidate how $\Lambda_\mathcal{L}^k$ in (2a) is determined based on (2b), let investigate the last term in (2c) (i.e., $(y_\mathcal{L}^{c,a} V_n^a + y_\mathcal{L}^{c,b} u_\mathcal{L}^b V_n^b + y_\mathcal{L}^{c,c} V_n^c)$). Based on (2b), if $k \neq \phi$ or $p$, $\Lambda_\mathcal{L}^k$ is equal to $u_\mathcal{L}^k$. In the last term of (2c), $\phi = a$ and $p = c$. Therefore, only $u_\mathcal{L}^b$ is appeared in this last term.

In (2a), if phase $\phi$ of the line is disconnected, $u_\mathcal{L}^\phi$ will be equal to zero, which results in $I_\mathcal{L}^\phi$ to be obtained as zero. Let assume phase $\phi$ of the line is connected (i.e., $u_\mathcal{L}^\phi = 1$). In this regard, if any of the other two phases of the line are disconnected, the binary variables corresponding to them (i.e., $u_\mathcal{L}^p$ or $u_\mathcal{L}^k$ in (2a)) will be equal to zero. The zero values for these binary variables in (2a) result in eliminating the effect of the voltages (i.e., $V_n^k$) and currents (i.e., $I_\mathcal{L}^p$) of the other two phases on the current of phase $\phi$ of the line (i.e., $I_\mathcal{L}^\phi$). For the given example in (2c), if phase $b$ of the line $\mathcal{L}$ is not energized, $u_\mathcal{L}^b$ will be obtained as zero. Since every $V_n^b$ and $I_\mathcal{L}^b$ is multiplied by $u_\mathcal{L}^b$ in (2c), the zero value of $u_\mathcal{L}^b$ results in vanishing of all terms including the voltage and current of phase $b$. Constraints (2d)-(2e) include the real and imaginary parts of $I_\mathcal{L}^\phi$ in (2a).

$I_\mathcal{L}^{r,\phi} = (R_\mathcal{L}^{\phi,\phi})^{-1} u_\mathcal{L}^\phi \left[ V_n^{r,\phi} - V_m^{r,\phi} + X_\mathcal{L}^{\phi,\phi} I_\mathcal{L}^{im,\phi} - \frac{1}{2} y_\mathcal{L}^{\phi,\phi} (R_\mathcal{L}^{\phi,\phi} V_n^{im,\phi} + X_\mathcal{L}^{\phi,\phi} V_n^{r,\phi}) - \frac{1}{2} (\sum_{k \in \psi, k \neq \phi} y_\mathcal{L}^{\phi,k} u_\mathcal{L}^k (R_\mathcal{L}^{\phi,\phi} V_n^{im,k} + X_\mathcal{L}^{\phi,\phi} V_n^{r,k})) - \sum_{p \in \psi, p \neq \phi} R_\mathcal{L}^{\phi,p} u_\mathcal{L}^p \left( I_\mathcal{L}^{r,p} + \frac{1}{2} (\sum_{k \in \psi} y_\mathcal{L}^{p,k} \Lambda_\mathcal{L}^k V_n^{im,k}) \right) + \sum_{p \in \psi, p \neq \phi} X_\mathcal{L}^{\phi,p} u_\mathcal{L}^p \left( I_\mathcal{L}^{im,p} - \frac{1}{2} \sum_{k \in \psi} y_\mathcal{L}^{p,k} \Lambda_\mathcal{L}^k V_n^{r,k} \right) \right], \forall \phi \in \psi, \mathcal{L} \in \mathcal{H}''$ (2d)

$I_\mathcal{L}^{im,\phi} = (R_\mathcal{L}^{\phi,\phi})^{-1} u_\mathcal{L}^\phi \left[ V_n^{im,\phi} - V_m^{im,\phi} - X_\mathcal{L}^{\phi,\phi} I_\mathcal{L}^{r,\phi} + \frac{1}{2} y_\mathcal{L}^{\phi,\phi} (R_\mathcal{L}^{\phi,\phi} V_n^{r,\phi} - X_\mathcal{L}^{\phi,\phi} V_n^{im,\phi}) + \frac{1}{2} (\sum_{k \in \psi, k \neq \phi} y_\mathcal{L}^{\phi,k} u_\mathcal{L}^k (R_\mathcal{L}^{\phi,\phi} V_n^{r,k} - X_\mathcal{L}^{\phi,\phi} V_n^{im,k})) - \sum_{p \in \psi, p \neq \phi} R_\mathcal{L}^{\phi,p} u_\mathcal{L}^p \left( I_\mathcal{L}^{im,p} - \frac{1}{2} (\sum_{k \in \psi} y_\mathcal{L}^{p,k} \Lambda_\mathcal{L}^k V_n^{r,k}) \right) - \sum_{p \in \psi, p \neq \phi} X_\mathcal{L}^{\phi,p} u_\mathcal{L}^p \left( I_\mathcal{L}^{r,p} + \frac{1}{2} \sum_{k \in \psi} y_\mathcal{L}^{p,k} \Lambda_\mathcal{L}^k V_n^{im,k} \right) \right], \forall \phi \in \psi, \mathcal{L} \in \mathcal{H}''$ (2e)

With the advent of DERs in modern distribution systems, more sensors such as micro-PMUs and smart meters are installed to improve the observability and control of the distribution systems. Let consider $\mathfrak{R}$ as a vector of measurements (e.g., micro-PMUs and smart meters) in an unbalanced distribution system, and $\mathfrak{R}(x)$ as a vector of their corresponding variables, which are functions of the system topology and states (i.e., $x$). The objective function of the proposed simultaneous topology processor and state estimation tool in an unbalanced distribution system is minimizing the weighted error between $\Re$ and $\Re(x)$ as expressed in (2f) using Euclidean norm (i.e., $\|\ \|_2$).

$$\bar{x}(\Re) = Arg\ min\ \|\mathcal{W}(\Re(x) - \Re)\|_2^2 \qquad (2f)$$

where $\bar{x}(\Re)$ is the vector of obtained status of switches and voltage phasors. $\mathcal{W}$ is a diagonal matrix including covariance of measurements noise. The objective function (2f) is subject to constraints (1d)-(1i), (2b), and (2d)-(2e). It should be noted that constraints (1d)-(1e) are considered for the lines without switch (i.e., $\mathcal{L} \in \mathcal{H}'$), and constraints (2d)-(2e) are modeled for lines equipped with the switch devices (i.e., $\mathcal{L} \in \mathcal{H}''$). Although the objective function in (2f) is convex; however, the constraints are nonlinear.

### C. Proposed Simultaneous Topology Processor and State Estimation Tool

The simultaneous topology processor and state estimation model in section II. B is a mixed-integer nonlinear quadratic programming (MINQP) problem, which may not lead to the global optimal solution. The nonlinearity in this problem is due to two sets of constraints (i.e., constraints (2d)-(2e) and constraints (1h)-(1i)). In order to attain the global optimal solution, these two sets of constraints are linearized in the following sections using the Big M method and the Tylor extension approach.

#### 1) Linearizing Line Current Constraints with Switches

The current of the distribution line $\mathcal{L} \in \mathcal{H}''$ with three single-phase switches is given in (2d) and (2e), in which the fifth, sixth, and seventh terms of each of them are nonlinear. The nonlinearity in the fifth term of (2d)-(2e) is due to the multiplication of binary variable $u_\mathcal{L}^k$ with continuous variables $V_n^{im,k}$ and $V_n^{r,k}$. Such nonlinearity is linearized by applying the big M method as shown in (3a)-(3d).

$$-M(1 - u_\mathcal{L}^k) \le d_k - V_n^{r,k} \le M(1 - u_\mathcal{L}^k) \qquad (3a)$$

$$-Mu_\mathcal{L}^k \le d_k \le Mu_\mathcal{L}^k \qquad (3b)$$

$$-M(1 - u_\mathcal{L}^k) \le d_k' - V_n^{im,k} \le M(1 - u_\mathcal{L}^k) \qquad (3c)$$

$$-Mu_\mathcal{L}^k \le d_k' \le Mu_\mathcal{L}^k \qquad (3d)$$

According to (3a)-(3d), if $u_\mathcal{L}^k = 1$, the auxiliary variables $d_k$ and $d_k'$ will be equal to $V_n^{r,k}$ and $V_n^{im,k}$, respectively. In the sixth and seventh terms of (2d)-(2e), if the axillary binary variable $\Lambda_\mathcal{L}^k$ is equal to 1 (i.e., when $k$ is equal to $\phi$ or $p$), multiplication of binary variable $u_\mathcal{L}^p$ with continuous variables $I_\mathcal{L}^{r,p}$, $I_\mathcal{L}^{im,p}$, $V_n^{r,k}$, and $V_n^{im,k}$ will be appeared in the constraints. To linearize these nonlinear terms, constraints (3e)-(3h) are added to the model.

$$-M(1 - u_\mathcal{L}^p) \le E_p - \left(I_\mathcal{L}^{r,p} + \frac{1}{2}\sum_{\substack{k \in \psi \\ k = \phi \vee p}} y_\mathcal{L}^{p,k} V_n^{im,k}\right) \le M(1 - u_\mathcal{L}^p) \qquad (3e)$$

$$-Mu_\mathcal{L}^p \le E_p \le Mu_\mathcal{L}^p \qquad (3f)$$

$$-M(1 - u_\mathcal{L}^p) \le E_p' - \left(I_\mathcal{L}^{im,p} - \frac{1}{2}\sum_{\substack{k \in \psi \\ k = \phi \vee p}} y_\mathcal{L}^{p,k} V_n^{r,k}\right) \le M(1 - u_\mathcal{L}^p) \qquad (3g)$$

$$-Mu_\mathcal{L}^p \le E_p' \le Mu_\mathcal{L}^p \qquad (3h)$$

where $E_p$ and $E_p'$ are auxiliary variables. In the linear constraints (3e)-(3h), if the switch installed at phase $p$ is disconnected (i.e., $u_{n,m}^p = 0$), $E_p$ and $E_p'$ will become zero. In this regard, the current of the line at phase $p$ and its corresponding mutual impedance (i.e., $R_{n,m}^{\phi,p}$ and $X_\mathcal{L}^{\phi,p}$) and admittances (i.e., $y_\mathcal{L}^{p,k}$) will be removed from the formulations of current at phase $\phi$ (i.e., (2d)-(2e)). If the axillary binary variable $\Lambda_\mathcal{L}^k$ is equal to $u_\mathcal{L}^k$ (i.e., when $k$ is not equal to $\phi$ or $p$) in the sixth and seventh terms of (2d)-(2e), multiplication of two binary variable $u_\mathcal{L}^p$ and $u_\mathcal{L}^k$ with continuous variables $V_n^{r,k}$ and $V_n^{im,k}$ will be imposed to the model. First, multiplication of two binary variables $u_\mathcal{L}^p$ and $u_\mathcal{L}^k$ is linearized by substituting it with a new binary variable $\mathscr{E}_\mathcal{L}^{p,k}$ and adding constraints (3i)-(3k).

$$\mathscr{E}_\mathcal{L}^{p,k} \le u_\mathcal{L}^p \qquad (3i)$$

$$\mathscr{E}_\mathcal{L}^{p,k} \le u_\mathcal{L}^k \qquad (3j)$$

$$\mathscr{E}_\mathcal{L}^{p,k} \ge u_\mathcal{L}^p + u_\mathcal{L}^k - 1 \qquad (3k)$$

Based on (3i)-(3k), if both $u_\mathcal{L}^p$ and $u_\mathcal{L}^k$ are equal to 1, $\mathscr{E}_\mathcal{L}^{p,k}$ will be obtained as 1. Let assume $u_\mathcal{L}^p = 0$ and $u_{n,m}^k = 0$, according to (3i)-(3j), $\mathscr{E}_\mathcal{L}^{p,k}$ is also forced to be zero.

Second, by replacing $u_\mathcal{L}^p \times u_\mathcal{L}^k$ with $\mathscr{E}_\mathcal{L}^{p,k}$, multiplication of $\mathscr{E}_\mathcal{L}^{p,k}$ with continuous variables $V_n^{im,k}$ and $V_n^{r,k}$ is linearized by introducing two auxiliary variables $H_k$ and $H_k'$ and the big M approach as defined in (3l)-(3o).

$$-M(1 - \mathscr{E}_\mathcal{L}^{p,k}) \le H_k - \frac{1}{2}y_\mathcal{L}^{p,k}V_n^{r,k} \le M(1 - \mathscr{E}_\mathcal{L}^{p,k}) \qquad (3l)$$

$$-M\mathscr{E}_\mathcal{L}^{p,k} \le H_k \le M\mathscr{E}_\mathcal{L}^{p,k} \qquad (3m)$$

$$-M(1 - \mathscr{E}_\mathcal{L}^{p,k}) \le H_k' - \frac{1}{2}y_\mathcal{L}^{p,k}V_n^{im,k} \le M(1 - \mathscr{E}_\mathcal{L}^{p,k}) \qquad (3n)$$

$$-M\mathscr{E}_\mathcal{L}^{p,k} \le H_k' \le M\mathscr{E}_\mathcal{L}^{p,k} \qquad (3o)$$

In the constraints (3l)-(3o), if $\mathscr{E}_\mathcal{L}^{p,k}$ is equal to 1, $H_k = \frac{1}{2}y_\mathcal{L}^{p,k}V_n^{r,k}$ and $H_k' = \frac{1}{2}y_\mathcal{L}^{p,k}V_n^{im,k}$. Considering the proposed linear constraints (3a)-(3o), the constraints (2d) and (2e) can be written as (3p) and (3.q), respectively.

$$I_\mathcal{L}^{r,\phi} = \left(R_\mathcal{L}^{\phi,\phi}\right)^{-1}u_\mathcal{L}^\phi\Big[V_n^{r,\phi} - V_m^{r,\phi} + X_\mathcal{L}^{\phi,\phi}I_\mathcal{L}^{im,\phi} - \frac{1}{2}y_\mathcal{L}^{\phi,\phi}\left(R_\mathcal{L}^{\phi,\phi}V_n^{im,\phi} + X_\mathcal{L}^{\phi,\phi}V_n^{r,\phi}\right) - \frac{1}{2}\left(\sum_{k\in\psi,k\ne\phi}y_\mathcal{L}^{\phi,k}\left(R_\mathcal{L}^{\phi,\phi}d_k' + X_\mathcal{L}^{\phi,\phi}d_k\right)\right) - \sum_{p\in\psi,p\ne\phi}R_\mathcal{L}^{\phi,p}\left(E_p + H_k'\right) + \sum_{p\in\psi,p\ne\phi}X_\mathcal{L}^{\phi,p}\left(E_p' - H_k\right)\Big], \forall\ \phi \in \psi, \mathcal{L} \in \mathcal{H}'' \qquad (3p)$$

$$I_\mathcal{L}^{im,\phi} = \left(R_\mathcal{L}^{\phi,\phi}\right)^{-1}u_\mathcal{L}^\phi\Big[V_n^{im,\phi} - V_m^{im,\phi} - X_\mathcal{L}^{\phi,\phi}I_\mathcal{L}^{r,\phi} + \frac{1}{2}y_\mathcal{L}^{\phi,\phi}\left(R_\mathcal{L}^{\phi,\phi}V_n^{r,\phi} - X_\mathcal{L}^{\phi,\phi}V_n^{im,\phi}\right) + \frac{1}{2}\left(\sum_{k\in\psi,k\ne\phi}y_\mathcal{L}^{\phi,k}\left(R_\mathcal{L}^{\phi,\phi}d_k - X_\mathcal{L}^{\phi,\phi}d_k'\right)\right) - \sum_{p\in\psi,p\ne\phi}R_\mathcal{L}^{\phi,p}\left(E_p' - H_k\right) - \sum_{p\in\psi,p\ne\phi}X_\mathcal{L}^{\phi,p}\left(E_p + H_k'\right)\Big], \forall\ \phi \in \psi, \mathcal{L} \in \mathcal{H}'' \qquad (3q)$$

In constraints (3p) and (3.q), all the binary variables except $u_\mathcal{L}^\phi$ are eliminated, but their impacts are modeled by including auxiliary variables $d_k$, $d_k'$, $E_p$, $E_p'$, $H_k$, and $H_k'$ based on constraints (3a)-(3o). Let consider inside the brackets in (3p) and (3q) to be $F_\phi^r$ and $F_\phi^{im}$. In this regard, constraints (care expressed as follows:

$$I_\mathcal{L}^{r,\phi} = \left(R_\mathcal{L}^{\phi,\phi}\right)^{-1} u_\mathcal{L}^\phi [F_\phi^r], \forall \phi \in \psi, \mathcal{L} \in \mathcal{H}'' \quad (3r)$$

$$I_\mathcal{L}^{im,\phi} = \left(R_\mathcal{L}^{\phi,\phi}\right)^{-1} u_\mathcal{L}^\phi [F_\phi^{im}], \forall \phi \in \psi, \mathcal{L} \in \mathcal{H}'' \quad (3s)$$

Although $F_\phi^r$ and $F_\phi^{im}$ are linear, $I_\mathcal{L}^{r,\phi}$ and $I_\mathcal{L}^{im,\phi}$ are still nonlinear because of the product of binary variable $u_\mathcal{L}^\phi$ with continuous auxiliary variables $F_\phi^r$ and $F_\phi^{im}$. To linearize the nonlinearity as a result of modeling the status of single-phase switch of the phase $\phi$ (i.e., $u_\mathcal{L}^\phi$) in (3r) and (3s), constraints (3t)-(3w) based on the big M method are considered.

$$-M(1-u_\mathcal{L}^\phi) \leq I_\mathcal{L}^{r,\phi} - \left(R_\mathcal{L}^{\phi,\phi}\right)^{-1}[F_\phi^r] \leq M(1-u_\mathcal{L}^\phi), \forall \phi \in \psi, \mathcal{L} \in \mathcal{H}'' \quad (3t)$$

$$-M(1-u_\mathcal{L}^\phi) \leq I_\mathcal{L}^{im,\phi} - \left(R_\mathcal{L}^{\phi,\phi}\right)^{-1}[F_\phi^{im}] \leq M(1-u_\mathcal{L}^\phi), \forall \phi \in \psi, \mathcal{L} \in \mathcal{H}'' \quad (3u)$$

$$-Mu_\mathcal{L}^\phi \leq I_\mathcal{L}^{r,\phi} \leq Mu_\mathcal{L}^\phi, \forall \phi \in \psi, \mathcal{L} \in \mathcal{H}'' \quad (3v)$$

$$-Mu_\mathcal{L}^\phi \leq I_\mathcal{L}^{im,\phi} \leq Mu_\mathcal{L}^\phi, \forall \phi \in \psi, \mathcal{L} \in \mathcal{H}'' \quad (3w)$$

Based on constraints (3t)-(3w), if $u_\mathcal{L}^\phi$ is equal to zero, $I_\mathcal{L}^{r,\phi}$ and $I_\mathcal{L}^{im,\phi}$ will become zero, which implies that phase $\phi$ of line between bus $n$ to bus $m$ is disconnected.

*2) Linearization of Power Balance Constraints*

The right side of power balance equations (1h)-(1i) is not linear due to the product of current and voltage variables. In this paper, iterative first-order approximation of Taylor extension is utilized in order to propose the linear model of (1h)-(1i) as follows:

$$\sum_{\forall g \in g(n)} P_{g,\phi}^G - \sum_{\forall l \in l(n)} d_{p,l,\phi} - \sum_{\forall r \in R(n)} NL_{r,\phi}^{TR} = V_{n,it-1}^{r,\phi} I_n^{r,\phi} + V_{n,it-1}^{im,\phi} I_n^{im,\phi} + I_{n,it-1}^{r,\phi} V_n^{r,\phi} + I_{n,it-1}^{im,\phi} V_n^{im,\phi} - V_{n,it-1}^{r,\phi} I_{n,it-1}^{r,\phi} - V_{n,it-1}^{im,\phi} I_{n,it-1}^{im,\phi}, \forall \phi \in \psi, n \in \mathcal{J} \quad (3x)$$

$$\sum_{\forall g \in g(n)} Q_{g,\phi}^G + \sum_{\forall c \in C(n)} Q_{c,\phi}^C - \sum_{\forall l \in l(n)} d_{q,l,\phi} = V_{n,it-1}^{im,\phi} I_n^{r,\phi} - V_{n,it-1}^{r,\phi} I_n^{im,\phi} + I_{n,it-1}^{r,\phi} V_n^{im,\phi} - I_{n,it-1}^{im,\phi} V_n^{r,\phi} - V_{n,it-1}^{im,\phi} I_{n,it-1}^{r,\phi} + V_{n,it-1}^{r,\phi} I_{n,it-1}^{im,\phi}, \forall \phi \in \psi, n \in \mathcal{J} \quad (3y)$$

where $V_{n,it-1}^{r,\phi}, V_{n,it-1}^{im,\phi}, I_{n,it-1}^{r,\phi}$ and $I_{n,it-1}^{im,\phi}$ are output of the previous iteration (i.e., $it-1$).

Finally, the proposed simultaneous state estimation and topology processor tool in an unbalanced distribution system with DERs and emerging sensors such as micro-PMUs is formulated as a MIQP problem using (3z).

$$\bar{x}(\Re) = arg\ min\ \|\mathcal{W}(\Re(x) - \Re)\|_2^2 \quad (3z)$$

*Subject to* $(1d) - (1g), (3a) - (3q), (3t) - (3y)$

It is worth to note that all constraints are linear and the objective function is convex. The structure of the proposed simultaneous topology processor and state estimation tool for an unbalanced distribution system in the real-time is illustrated in Fig. 1. It should be noted that the voltage measurements are obtained by solving an AC power flow via OpenDSS [20]. As shown in Fig. 1, the measurements with added noise are fed to the proposed tool. In the first iteration, the proposed tool uses the flat start point. The more accurate values of $V_{n,it-1}^{r,\phi}, V_{n,it-1}^{im,\phi}, I_{n,it-1}^{r,\phi}$ and $I_{n,it-1}^{im,\phi}$ can be obtained by solving the proposed simultaneous topology processor and state estimation tool iteratively, which results in improving the linearization approximations in the model.

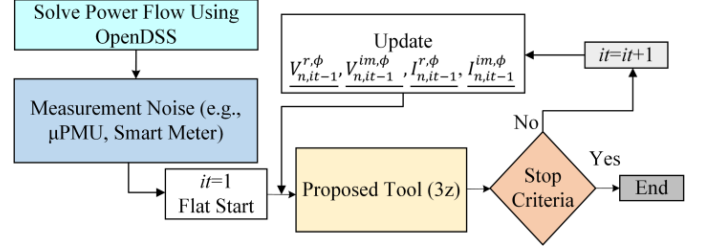

Fig. 1. Proposed topology processor and state estimation in an unbalanced distribution system.

### III. ANALYTIC FORMULATION OF REAL-RIME SIMULTANEOUS STATE ESTIMATION AND TOPOLOGY AND OUTAGE DETECTION

Identifying islands and outages in an unbalanced distribution system is challenging due to the limited observability of the system. The status of disconnected lines, which lead to the outage and create the islands in the system, must be zero. However, in (3v)-(3w), when $I_\mathcal{L}^{r,\phi} = 0$ and $I_\mathcal{L}^{im,\phi} = 0$, the status of switch at phase $\phi$ of the line ($u_\mathcal{L}^\phi$) can be identified as 0 or 1. Hence, the proposed topology processor and state estimation model given in (3z) is improved by adding following linear constraints in order to identify outages in an unbalanced distribution system.

$$\vartheta - I_\mathcal{L}^{r,\phi} \leq \lambda (1 - \alpha_\mathcal{L}^\phi) \quad (4a)$$

$$\vartheta + I_\mathcal{L}^{r,\phi} \leq \lambda (1 - \beta_\mathcal{L}^\phi) \quad (4b)$$

$$\vartheta - I_n^{im,\phi} \leq \lambda (1 - \gamma_\mathcal{L}^\phi) \quad (4c)$$

$$\vartheta + I_n^{im,\phi} \leq \lambda (1 - \mu_\mathcal{L}^\phi) \quad (4d)$$

$$u_\mathcal{L}^\phi \leq \alpha_\mathcal{L}^\phi + \beta_\mathcal{L}^\phi + \gamma_\mathcal{L}^\phi + \mu_\mathcal{L}^\phi \quad (4e)$$

where $\alpha_\mathcal{L}^\phi, \beta_\mathcal{L}^\phi, \gamma_\mathcal{L}^\phi$, and $\mu_\mathcal{L}^\phi$ are auxiliary binary variables, $\vartheta$ is a small positive number, and $\lambda$ is a large positive number. According to (4a)-(4e), if $I_\mathcal{L}^{r,\phi}$ and $I_\mathcal{L}^{im,\phi}$ are zero, $u_\mathcal{L}^\phi$ will be forced to be zero. The formulation of outage, topology, and state identification model is defined in (4f).

$$\bar{x}(\Re) = arg\ min\ \|\mathcal{W}(\Re(x) - \Re)\|_2^2 \quad (4f)$$

*Subject to* $(1d)-(1g), (3a)-(3q), (3t)-(3y), (4a)-(4e)$

### IV. SENSORS DATA INTEGRATION APPROACH IN DISTRIBUTION SYSTEM

Smart meters typically gauge electric energy consumptions. They can also provide average active power but not reactive power. To address lack of reactive power measurement in the distribution system, pseudo reactive power measurements based on the random power factor are considered in the literature [7]. However, pseudo reactive power measurements introduce high errors in the model. To overcome the limitation of smart meters not providing reactive power, constraint (5a) is embedded in the proposed topology, outage, and state identification model.

$$P_i \sqrt{\left(\frac{1}{PF_{i,max}^2} - 1\right)} \leq Q_i \leq P_i \sqrt{\left(\frac{1}{PF_{i,min}^2} - 1\right)}, \forall i \in S \quad (5a)$$

where $P_i$ is the variable associated with the average active power measured by the smart meter $i$, $Q_i$ is reactive power demand at each load, and $S$ is the set of smart meters. $PF_{i,min}$ and $PF_{i,max}$ is the minimum and maximum power factor of each load and are considered 0.9 and 0.99, respectively.

Another challenge in the distribution system studies is that micro-PMUs provide synchronized instantaneous measurements of voltage and current phasors with the frequency of 120Hz [21]; however, smart meters can provide an average of energy consumptions over a time period of 15 minutes, which are not instantaneous measurements. To address the challenge of integrating measurements of these two devices with different sampling rate and approaches (i.e., average versus instantaneous), a statistical study is performed in this paper. In this regard, the residential load data of Pecan street Inc. database are considered [22]. Since the proposed model is conducted for the primary distribution system, smart meter data are aggregated at the secondary side of each distribution transformer. Therefore, 1 second load data of different groups of customers (e.g., 7-8 houses) are collected at the secondary side of distribution transformers. Based on the conducted statistical study, the load variation between average active power during a 15 minutes time interval and instantaneous active power measurements follows a Gaussian distribution $\mathcal{N} = (0, \sigma^2)$ with zero mean and standard deviation $\sigma$, where 99.73% of data are set in $3\sigma = 30\%$. Therefore, 30% error is considered due to utilizing average power measured by smart meters versus instantaneous active power measurements due to the device sampling limitations.

## V. SIMULATION RESULTS

The proposed model is simulated using the distribution feeder of a local electric utility in Arizona. The primary feeder of the test system as shown in Fig. 2 contains 2100 buses, 2161 distribution lines, and 371 distribution transformers. There are 23 switch cabinets with 157 single-phase switches inside. The total number of single-phase switches is 859, including 681 single-phase transformer switches, two three-phase transformer switches, four three-phase capacitor bank switches, and one three-phase switch at the substation. The number of aggregated loads with smart meters at the secondary of distribution transformers is 342, and the number of aggregated PV with smart meter data is 251. Measurement data of the distribution feeder of the local electric utility include one pseudo micro-PMU measurement, smart meter data of load buses, rooftop PV output power, substation measurements, and reactive power injection of capacitors. The micro-PMU is placed in the substation bus. The Gaussian distribution function is utilized to model the measurement noise of the pseudo micro-PMU with $Total\ Vector\ Error \leq 0.05\%$ [6]. The measurement data of smart meters, rooftop PV sensors, capacitors, and substation sensor device include noise, which is considered as a Gaussian distribution with zero mean and 10%, 10%, 1%, and 1% error, respectively [23]. Topology processor in utility distribution feeders is more challenging compared to IEEE test case systems due to switch cabinets with different switch configurations, as shown in Fig. 3, various phase configurations (i.e., single-phase, two-phase, and three-phase) of feeders and laterals as demonstrated in Fig. 2, and other detailed models of the electric devices. To describe the test system in detail, switch cabinets with various phase configurations and no-load loss of transformers are modeled in the proposed simultaneous topology processor and state estimation tool.

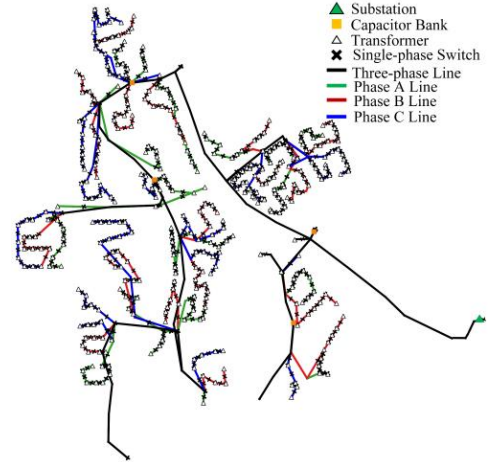

Fig. 2. Distribution feeder of a local electric utility in Arizona.

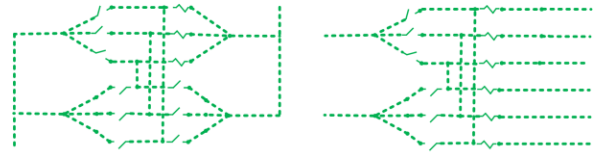

Fig. 3. (Left) switch cabinet with 9 single-phase switches and three fuses, (Right) switch cabinet with 6 single-phase switches and 6 fuses.

Table I shows the accuracy of the proposed simultaneous topology processor and state estimation tool in identifying the status of switches of the test system in a normal state of the system without outage. Two metrics are used in Table I to assess the accuracy of the proposed tool in detection of of the topology of the system. The first metric is calculated using $\left(M_1 = 100 \times \left(1 - \left[\sum_{\varepsilon \in E'} \sum_\phi \left|u_\varepsilon^\phi - u_\varepsilon^{\phi'}\right|/S'\right]\right)\right)$, and represents the accuracy of the proposed tool without considering the status of switches connected to presently no-load region (i.e., not connected to any device in the system) and normally open switches (NOSs). Set $E'$ includes $\{1,..,S'\}$ in which $S'$ is the total number of single-phase switches excluding those connected to presently no-load region and NOSs. The second metric is calculated using $\left(M_2 = 100 \times \left(1 - \right.\right.$

$\left[\sum_{\tau \in E} \sum_{\phi} \left|u_\tau^\phi - u_\tau^{\phi'}\right|/S\right]\right)$, where set $E = \{1,..,S\}$ and $S$ is the total number of single-phase switches in the system. $M_2$ shows the accuracy of the proposed model considering the status of all switches in the system even those connected to pr-sently no-load region and NOSs. Location of 12 errors reported in Table I is demonstrated in Fig. 4 using five orange circles. Yellow rhombuses in Fig. 4 show ending of the laterals connected to presently no-load region. It can be observed that most of the errors (i.e., 9 errors in circles 3-5) are inside switch cabinets that are located in parts of the system, which are connected to presently no-load regions. In other words, they are not connected to any device (e.g., load, distribution transformer, capacitor, rooftop PV unit) in the network. These locations are related to where the local electric utility has extended the feeder to serve future customers, but no development has taken place in those regions yet. The remaining three errors are shown in circles 1-2, which are in parts of the system that a NOS connects different phases together. While the proposed algorithm identifies the status of NOSs correctly as disconnected, the status of switches right before the NOSs are identified wrongly as disconnected due to current of zero. It can be observed in Table I that the accuracy of the proposed model is high for the primary topology processor with few numbers of measurement devices (e.g., one micro-PMU) and consideration of measurement noise with and without considering the status of switches connected to presently no-load region and NOSs.

TABLE I. ACCURACY OF TOPOLOGY PROCESSOR WITHOUT OUTAGE.

| Measurement noise | Yes |
|---|---|
| Total number of single-phase switches | 859 |
| Wrong statuses (connected to presently no-load region and NOS) | 12 |
| Undetectable statuses | 0 |
| $M_1$ accuracy (%) | 100 |
| $M_2$ accuracy (%) | 98.60 |
| Captured feeder load (%) | 100 |
| Simulation time (s) | 20 |

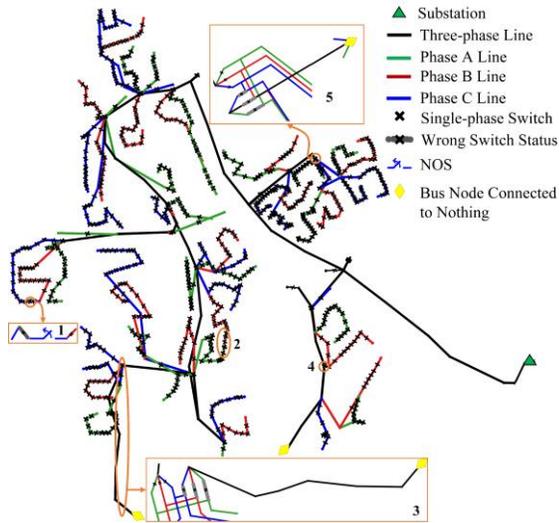

Fig. 4. Location of 12 topology processor errors.

With high interconnection structure accuracy (i.e., 100% or 98.6%), 100% of local electric utility feeder load is captured as shown in Table I. The state estimation results including voltage magnitude error and voltage angle error for each node of the local electric utility network are depicted in Fig. 5. In Fig. 5, each point represents the voltage error of phase a, b, or c of buses. According to Fig. 5, the proposed simultaneous topology processor and state estimation tool estimates states of the test system very precisely during a normal condition (maximum voltage magnitude error $< 0.2\%$ and maximum voltage angle error $< 0.002$ radian). The average simulation time is 20 seconds, which implies that the proposed tool is fast for the real-time simulations. In order to further illustrate the performance of the proposed simultaneous topology processor and state estimation tool, four other sets of cases are investigated in the following. These four case studies include studying outages in the test system, modeling bad data, practical consideration of micro-PMU and smart meter devices, and evaluating the robustness of the proposed tool against the missing data.

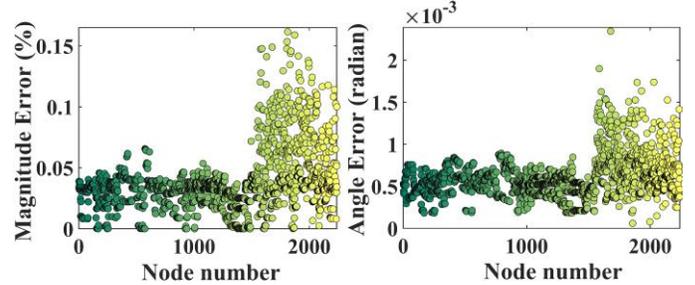

Fig. 5. State estimation error for a normal state of the local electric utility system (Left) voltage magnitude (Right) voltage angle for each node.

*A. Simultaneous Identification of Outage, Topology, and States of Distribution Systems*

The outage in the distribution system can occur due to switch malfunctions or faults in the system. As a result of the outages in a distribution system, the network is divided into the different energized regions and de-energized parts. The proposed model is able to detect not only the outage areas but also the topology of the energized regions. Four cases are analyzed for simultaneous detecting outage areas and topology of energized sections in the test system. In cases A1 and A2, one switch inside two different switch cabinets in the test system feeder is disconnected leading to two different outage areas in the system with 46 and 37 islanded nodes, respectively. In case A3, multiple simultaneous outages are simulated by considering cases A1 and A2 together. In case 4, a large outage as a result of disconnecting a switch inside another switch cabinet in the test system is simulated with 381 islanded nodes as illustrated in Fig. 6 (the red color shows the parts of the system with outage). Table II shows the accuracy of the proposed topology processor model for four simulated outage cases. As shown in Table II, the proposed model can detect topology and outages of the test system with a high accuracy in cases A1-A4. It should be noted that the locations of 12 errors in the status of single-phase switches in Table II for the four cases are the same as locations shown in Fig. 4. Moreover, the undetectable switch status in cases A1 and A3-A4 in Table II is the switch located right before the disconnected switch in each case. While this switch is connected in each case, no current is flowing on it due to the disconnection of the downstream feeder. The current of zero makes the status of such switches undetectable. However, the status of this switch does not impact load and feeder connectivity and the ability to control the feeder. With the high accuracy of the proposed approach in detecting the outage and topology of the system, 100% of the local electric utility feeder load is captured in cases A1-A4 for the energized areas. Due to the fact that the proposed tool identifies the system states simultaneously with the topology and outage areas of the system, the error of obtained voltage magnitude and angle states are illustrated in Fig. 7 for case A4 with the largest outage. It can be seen that the system states are estimated with the high accuracy even while having a large outage in the network.

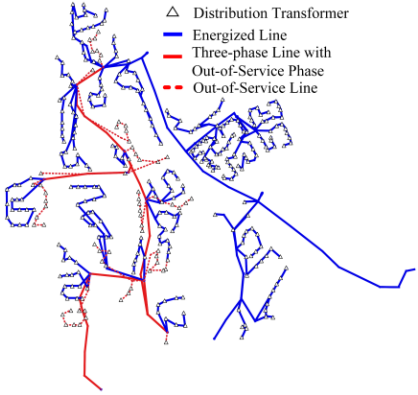

Fig. 6. Islanded region of the local electric utility system in case A4.

TABLE II. RESULTS OF TOPOLOGY PROCESSOR FOR FOUR CASES WITH OUTAGE.

| Case | A1 | A2 | A3 | A4 |
|---|---|---|---|---|
| Number of islanded nodes | 46 | 37 | 83 | 381 |
| Measurement noise | Yes | Yes | Yes | Yes |
| Total number of single-phase switches | 859 | 859 | 859 | 859 |
| Wrong statuses (connected to presently no-load region) and NOS | 12 | 12 | 12 | 12 |
| Undetectable statuses | 1 | 0 | 1 | 1 |
| $M_1$ accuracy (%) | 99.88 | 100 | 99.88 | 99.88 |
| $M_2$ accuracy (%) | 98.49 | 98.60 | 98.49 | 98.49 |
| Captured feeder load (%) | 100 | 100 | 100 | 100 |

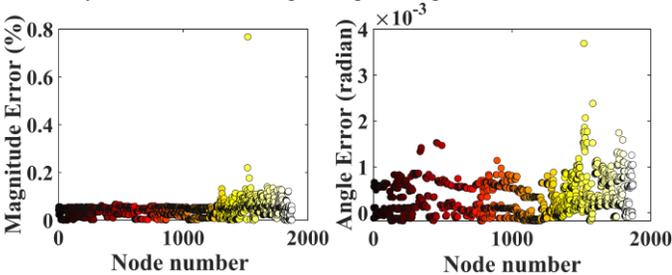

Fig. 7. State estimation error for case A4 with the largest outage (Left) voltage magnitude (Right) voltage angle for energized nodes of the system.

### B. Bad Data Modeling

In order to evaluate the robustness of the proposed tool in the case of bad data, three cases are studied, where the aggregated smart meter data at the secondary of the distribution transformers are vitiated for the bad data modeling. In this regard, in addition to the measurement noise of devices in the system (e.g., micro-PMU noise, 10% error in smart meter data), 60% error in the aggregated smart meter measurements connected to the largest loads in the test system is modeled for the bad data modeling. For instance, the active power error in all aggregated loads of the system in case B3 is shown in Fig. 8. As shown in Fig. 8, 40 aggregated smart meters data have an error equals to 60%, while the rest of the aggregated smart meters data are modeled with a 10% error due to the noise of the device using the Gaussian distribution. Table III shows the results of the proposed topology processor in the case of bad data in the system. As illustrated in Table III, both $M_1$ and $M_2$ metrics for evaluation of the topology processor accuracy in cases B1-B3 are the same as results of table I without considering bad data in the measurements. Therefore, it can be inferred that the proposed tool is significantly robust in performing topology processor in the case of considering both measurement noise and bad data. At the same time, the proposed approach with the high topology detection accuracy captures 100% of the local electric utility feeder load in cases B1-B3, as shown in Table III.

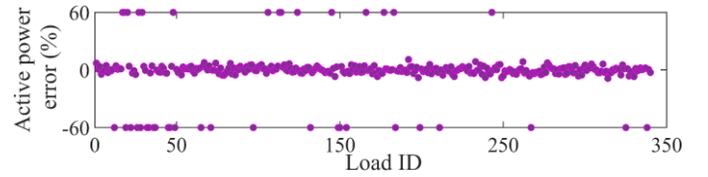

Fig. 8. Modeled smart meters measurements error for case B3 considering device error and bad data.

TABLE III. RESULTS OF THE PROPOSED TOPOLOGY PROCESSOR FOR THREE CASES WITH BAD DATA MODELING.

| Case | B1 | B2 | B3 |
|---|---|---|---|
| Number of aggregated smart meter measurements modeled as bad data | 10 | 30 | 40 |
| Measurement noise | Yes | Yes | Yes |
| Total number of single-phase switches | 859 | 859 | 859 |
| Wrong statuses (connected to presently no-load region) and NOS | 12 | 12 | 12 |
| Undetectable statuses | 0 | 0 | 0 |
| $M_1$ accuracy (%) | 100 | 100 | 100 |
| $M_2$ accuracy (%) | 98.60 | 98.60 | 98.60 |
| Captured feeder load (%) | 100 | 100 | 100 |

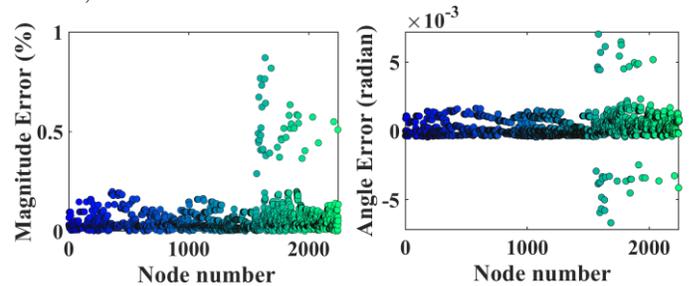

Fig. 9. State estimation error for case B3 (Left) voltage magnitude (Right) voltage angle for each node of the local electric utility network.

The error of estimated values of system states for case B3 with the highest number of aggregated smart meter measurements with the bad data (i.e., 40) are depicted in Fig. 9. By comparing Fig. 5 with Fig. 9, it can be observed that state estimation errors for most of the nodes in the test system are comparable, and only the error for those nodes with the bad data is slightly higher (e.g., maximum voltage magnitude error in case B3 < 0.9%). This comparison implies that the proposed tool is also robust in terms of estimating system states considering bad data in the system measurements.

## C. Practical Challenges of Utilizing Smart Meter and Micro-PMU in Distribution Systems

In the results of Table I-III, it is assumed that the smart meters can provide reactive power at the same time that the micro-PMU measurements are received at the feeder head. In this section, the practical challenges of exploiting the smart meters and the micro-PMUs measurements in the distribution system are examined: (1) smart meters usually measure the electric energy consumption over a 15 minutes interval without the reactive power of load; (2) integration of instantaneous synchronized micro-PMUs measurements with unsynchronized average smart meters measurements. In order to address these challenges, constraint (5a) is added to the proposed model for estimating the reactive power of the smart meters. Moreover, in addition to the 10% measurement error of the smart meter devices, 30% error using the Gaussian distribution is modeled for the smart meters measurements due to integrating them with the micro-PMU measurements with different sampling rates and approaches. The modeled active load deviations of the smart meters of the test system are shown in Fig. 10. It should be noted that all measurement noise of the devices in the system (e.g., micro-PMU, PV sensors) are also modeled.

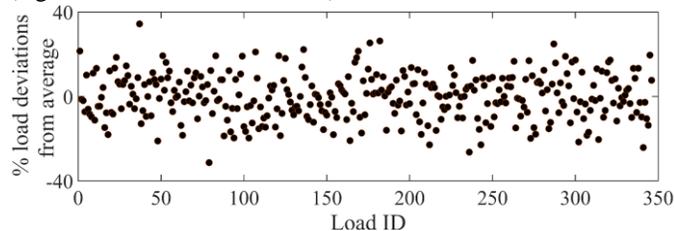
Fig. 10. Active load deviations from 15-min average smart meter measurements while considering device measurement error.

In this case, the total number of undetectable and wrongly identified (connected to presently no-load region and NOSs) status of switches are zero and 12, respectively. Accordingly, $M_1$ is obtained as 100% and $M_2$ is equal to 98.60%. The 12 wrong statues, in this case, are the same as those presented in Table I without considering the practical challenges of smart meters and micro-PMUs. These same results convey that the proposed primary topology processor model is highly precise even with considering measurements noise of all sensors and integration of synchronized instantaneous measurements of micro-PMU with average active power measurements of smart meter devices. Also, the proposed approach captures 100% of the system load. The error of state estimation using the proposed model for each node of the system is demonstrated in Fig. 11. Comparing Fig. 5 (i.e., results of the case without considering practical challenges of smart meters and micro-PMUs) with Fig. 11 illustrates that the accuracy of the proposed simultaneous topology processor and state estimation tool even while considering measurement devices with different sampling rates and approaches is significantly high.

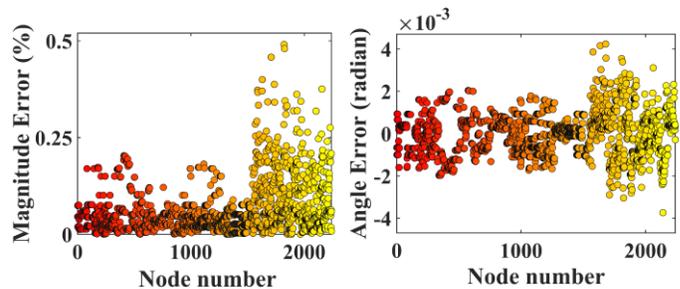
Fig. 11. Error for (Left) voltage magnitude and (Right) voltage angle estimation considering practical challenges of smart meters and micro-PMU.

## D. Missing Data Analysis

Another practical challenge of the real-time topology processor and state estimation tool is missing AMI data due to communication issues or the smart meter malfunction. There are around 6-8 smart meters aggregated at the secondary of the distribution transformers in the test system. In order to analyze the robustness of the proposed real-time tool, a case is considered, where some of the smart meter reading data at 246 buses (out of 342 load buses) are missing randomly. It is assumed that the maximum missed smart meter reading data is 50%. In this case, all measurement noise of the sensor devices (e.g., micro-PMU, smart meter) in the test system are modeled. Also, all other practical challenges presented in section III-C including 30% error for integration of the smart meters with the micro-PMU and inability of the smart meters for providing the reactive power of the loads are considered.

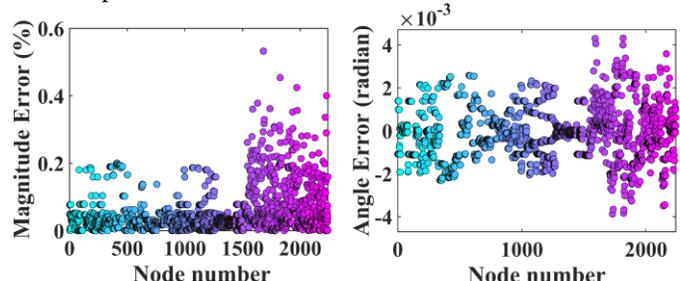
Fig. 12. Error for (Left) voltage magnitude and (Right) voltage angle estimation considering missed smart meter data.

In this paper, an approach is proposed for dealing with the missing data to enhance the accuracy of the proposed tool. In this regard, the missed smart meter measurement of each load is estimated based on the measurement of the previous time interval. It is assumed that the load deviation between two consecutive time intervals is maximum 30%. Therefore, a 30% load deviation using a Gaussian distribution with zero mean is considered for measurements of the previous time interval to estimate the reading data of the smart meter in the current time interval. The results of the proposed topology processor considering estimating missed smart meter data approach are $M_1 = 100\%$ and $M_2 = 98.60\%$. The total number of undetectable statues is zero. The total number of wrong statues is 12, which their locations are shown in Fig. 4. Even with considering about 72% (i.e., 246 out of 342) of load buses with the missed smart mere reading data, 100% of the local utility feeder load is captured using the proposed estimating missed smart meter data approach. The error of estimated states of the local electric utility system are depicted in Fig. 12. Comparing

Fig. 5, Fig. 11, and Fig. 12 confirm that the proposed approach is robust to missing data with considering all measurement noise of the sensors and the practical challenges of the smart meters and the micro-PMU.

## VI. Conclusion

In this paper, an efficient MIQP-based optimization tool is proposed to simultaneously identify network topology, estimate system state and detect outages of the unbalanced distribution systems. An ACOPF approach based on current and voltage is developed, which models the distribution system precisely. Moreover, the challenges of integrating micro-PMU and smart meter data with different sampling rates and approaches (i.e., average versus instantaneous) is studied using a statistical study. The proposed tool is further developed to overcome the limitation of the lack of reactive power measurements from smart meters. The results illustrate that the proposed model can identify different topologies and outages of the local electric utility system in Arizona with a high accuracy greater than 98% considering the challenges and measurements noise of the sensors in the system. Also, the proposed tool identifies the states of the system simultaneous with topologies and outages in the real-time accurately and fast (average simulation time is 20 seconds). The robustness of the proposed approach is evaluated by modeling bad data in the sensors' measurements. The results show that the proposed algorithm is significantly robust against bad data in the system measurements in terms of identifying the topology and states of the system. Also, 100% of the local utility feeder load is captured using the proposed tool in all case studies.